\documentclass[a4paper]{article}
\usepackage[latin1]{inputenc}
\sf
\usepackage{amssymb}
\def\R{\mathbb{R}}

\def\C{\mathbb{C}}

\def\Z{\mathbb{Z}}

\def\pr{$\bf{Proof :}$}

\def\vp{\varphi}

\def\R{\mathbb{R}}

\def\C{\mathbb{C}}

\def\Z{\mathbb{Z}}

\def\la{\lambda}

\def\Vect{{\rm Vect}\,}

\newcommand{\fila}[2]{\mbox{\scriptsize \vbox{\hbox{$#1$}\hbox{$#2$}}}}

\newtheorem{theo}{Theorem}[section]

\newtheorem{cor}[theo]{Corollary}
\newtheorem{lem}[theo]{Lemma}

\begin{document}

\title{On the Leibniz cohomology of vector fields}     

\author{Alessandra Frabetti\thanks{
I.R.M.A. - Universit\'e L. Pasteur, 
7 rue Ren\'e Descartes,
F-67084 Strasbourg Cedex. 
{\tt frabetti@math.u-strasbg.fr}}, 
Friedrich Wagemann\thanks{
I.G.D. - Universit\'e Lyon-I,  
43, bd du 11 Novembre 1918 
F-69622 Villeurbanne Cedex. 
{\tt wagemann@desargues.univ-lyon1.fr}}} 

\date{May 15, 2000}

\maketitle


\begin{abstract}
I.~M.~Gelfand and D.~B.~Fuks have studied the cohomology of the Lie algebra 
of vector fields on a manifold. 
In this article, we generalize their main tools to compute 
the Leibniz cohomology, by extending the two spectral sequences 
associated to the diagonal and the order filtration. 
In particular, we determine some new generators for the diagonal 
Leibniz cohomology of the Lie algebra of vector fields on the circle.
\end{abstract}


\section{Introduction}

Let $g$ be a Lie algebra. The Leibniz cohomology $HL^{*}(g)$ 
of $g$ with trivial coefficients is the cohomology of the complex with 
cochain modules $CL^{n}(g)={\rm Hom}(g^{\otimes n},\R)$ and coboundary 
operators 
$d: CL^{n}(g) \to CL^{n+1}(g)$, 
defined by 
$$
df(x_{1},\ldots,x_{n+1})= 
\sum_{i<j} (-1)^{j+1} f(x_{1},\ldots,x_{i-1},[x_{i},x_{j}],x_{i+1},\ldots,
\hat{x}_{j},\ldots,x_{n+1}),  
$$
where $(x_{1},\ldots,x_{k})$ denotes the element 
$x_{1} \otimes \ldots \otimes x_{k}$ in the tensor product 
$g^{\otimes k}$. The Leibniz cohomology 
is a generalization of the Chevalley-Eilenberg cohomology of Lie algebras, 
in the sense that it constitutes the natural cohomology theory for 
Leibniz algebras, which are 
non-antisymmetric generalizations of Lie algebras. 
It was introduced by J.-L. Loday and T. Pirashvili in \cite{LP}. 
When computed on Lie algebras, the Leibniz cohomology may detect new 
invariants in dimensions higher then $1$, since 
$HL^{0}(g)=H^{0}(g)\cong \R$ and 
$HL^{1}(g)=H^{1}(g)$. 
The method to compare the Lie and the Leibniz cohomology of a Lie 
algebra in higher dimensions is given by the Pirashvili long exact sequence 
$$
\begin{array}{l}
	0 \longrightarrow H^2(g) \longrightarrow HL^2(g) 
	\longrightarrow H_{rel}^0(g) \longrightarrow H^{3}(g) 
	\longrightarrow \ldots \\ 
	\hskip1cm \ldots \longrightarrow H^n(g) \longrightarrow HL^n(g) 
	\longrightarrow H_{rel}^{n-2}(g) \longrightarrow H^{n+1}(g) 
	\longrightarrow \ldots 
\end{array}
$$
induced by the projections $g^{\otimes n} \longrightarrow \Lambda^{n}g$, 
where $H_{rel}^*(g)$ is a relative cohomology described in \cite{P}. 

Loday and J.-M.~Oudom showed in \cite{L}, \cite{O} that the Leibniz cohomology 
$HL^{*}(g)$ is endowed with the structure of a graded dual Leibniz algebra 
(in the sense of Koszul duality for operads). 
In fact, there is a cup product 
$\cup: HL^{p}(g) \otimes HL^{q}(g) \longrightarrow HL^{p+q}(g)$ 
defined on the cohomology classes of two cocycles $\alpha \in CL^{p}(g)$, 
$\beta \in CL^{q}(g)$ as the cohomology class of the cocycle 
$$
\alpha \cup \beta(x_{1},\ldots,x_{p+q})= \!\!\!\!\!
	\sum_{\sigma \in Sh_{p-1,q}} \!\!\!\!\! {\rm sgn}(\sigma) 
	\alpha(x_{1},x_{\sigma (2)},\ldots,x_{\sigma (p)}) 
	\beta(x_{\sigma(p+1)},\ldots,x_{\sigma(p+q)}),   
$$
where $Sh_{p-1,q}$ denotes the set of the $(p-1,q)$-shuffles among 
the permutations on $\{2,3,...,p+q\}$. 
The cup product satisfies the relation of dual Leibniz algebras 
$$
	(\alpha \cup \beta) \cup \gamma - \alpha \cup (\beta \cup \gamma) 
	= (-1)^{|\beta| |\gamma|} \alpha \cup (\gamma \cup \beta) . 
$$ 
Therefore, in order to describe the Leibniz cohomology, it suffices 
to give its generators as a dual Leibniz algebra. 

When $g$ is a topological Lie algebra, we can consider the submodules of 
continuous cochains
$CL^{n}_{cont}(g) ={\rm Hom}_{cont}(g^{\otimes n},\R) \subset 
{\rm Hom}(g^{\otimes n},\R)=CL^{n}(g)$. 
Since the Leibniz differential $df$ of any continuous cochain $f$ is still
continuous, the sequence of continuous cochains 
$d:{\rm Hom}_{cont}(g^{\otimes n},\R) \to {\rm Hom}_{cont}(g^{\otimes
n+1},\R)$ forms a subcomplex of the Leibniz complex, whose cohomology 
is called continuous Leibniz cohomology of $g$. It is still denoted 
by $HL^{*}(g)$ when no confusion can arise. 

If $g$ is the topological Lie algebra $\Vect M$ of vector fields on
a differentiable manifold $M$, the Leibniz cohomology $HL^{*}(\Vect M)$ 
is a generalization of the Lie cohomology $H^{*}(\Vect M)$ studied by 
I.~M.~Gelfand, D.~B.~Fuks, R.~Bott, G.~Segal and A.~Haefliger in 
\cite{GelFuk}, \cite{Fuk}, \cite{BotSeg}, \cite{Hae} and references therein. 

In order to study the Leibniz cohomology of $\Vect \R^{n}$, 
J.~Lodder starts with the Lie algebra of formal vector fields 
$W_{n}=\R[[x_{1},...,x_{n}]]
\otimes \{\frac{\partial}{\partial x_{1}},..., \frac{\partial}
{\partial x_{n}}\}$, where $\R[[x_{1},...,x_{n}]]$ denotes 
the algebra of formal power series in $n$ variables, and 
$\frac{\partial}{\partial x_{i}}$ the partial derivation in respect of the 
variable $x_{i}$. 
He showed in \cite{Lod} that the first non-zero class in $HL^*(W_n)$ is the
Godbillon-Vey class in dimension $2n+1$. 
Using Pirashvili's sequence, he describes $HL^{*}(W_{n})$ as a dual 
Leibniz algebra with generators and relations given in \cite{Lod}. 
For $n=1$, this dual Leibniz algebra has one generator in degree $3$ 
(the Godbillon-Vey class $\theta$) and a new generator in degree $4$ 
(called $\alpha$) with relations $\theta^2=0$ and
$\alpha \cup \theta=0$. 

The case of $HL^{*}(\Vect M)$, for general smooth manifolds, is much
more complicated. Even in the simplest case $M=S^{1}$ there are too
many non-zero terms appearing in the sequence, which prevent us to 
determine the dimension of the Leibniz cohomology groups. 
Hence Pirashvili's method cannot be reproduced in general. 

In this paper, instead, we generalize the classical method of Gelfand
and Fuks, regarding continuous Leibniz cochains on $\Vect M$ as 
generalized sections of the tensor powers of the tangent bundle $TM$, over 
the manifold $M \times...\times M$. This guarantees for
Leibniz cohomology the same general results previously determined by 
Gelfand and Fuks for Lie cohomology, and some explicit computations 
for manifolds with few non-trivial de Rham cohomology classes. 

In particular, in section 2 we describe the spectral sequence
associated to the diagonal filtration (filtration of a
complex of distributions by their {\it support}). The first term of
the sequence is given by some quotient complexes.

In section 3 we describe a spectral sequence
converging to the cohomology of these quotient complexes. It is
associated to the order filtration, i.e. the filtration of a complex
of distributions which are supported on a given submanifold by their
{\it order}.

In the last section, we show that the order spectral sequence for the 
Leibniz cohomology of $\Vect S^1$ collapses. The resulting
diagonal cohomology is $4$-periodic, with dimension respectively 
$0, 1, 2, 1$. Finally, we describe the diagonal cocycles which determine 
the new generators.

For the collapsing of the diagonal spectral sequence, and thus for an
explicit desciption of all continuous Leibniz cohomology of
$\Vect S^1$, we would need the dual Leibniz algebra structure which is
still a mystery for us. The knowledge of this structure leads to 
study the relationship with Leibniz minimal models (as introduced by 
M.~Livernet in \cite{Mur}) and `real' homotopy types in the spirit of 
Bott, Haefliger and Segal \cite{Hae} \cite{BotSeg}. 


\section{Generalized sections and diagonal filtration}

Let $M$ be a smooth oriented manifold of dimension $n$.
Let $\Vect M$ be the Lie algebra of vector fields on $M$, equipped
with the $C^{\infty}$ topology which makes it a topological 
Fr\'echet nuclear Lie algebra. 
Let $CL^*_{cont}(\Vect M)$ be the complex of continuous
Leibniz cochains of $\Vect M$. 
Then, for any $m \geq 0$, the module of Leibniz $m$-cochains 
is the space of generalized sections of the bundle of tensor powers
$TM^{\otimes m}$ over the cartesian product manifold $M^m$, meaning
$TM^{\otimes m}=\bigotimes_{i=1}^m\pi^*_i(TM)$ where $\pi_i:M^m\to M$ is
the projection on the $i$th factor. 
In other words, the Leibniz $m$-cochains of $\Vect M$ are the 
distributions on $M^m$ with values in the bundle $TM^{\otimes m}$. 
Remark that the Chevalley-Eilenberg $m$-cochains of $\Vect M$, that is, 
the elements of $C^m_{cont}(\Vect M)$, can be seen as conveniently 
anti-symmetrized generalized sections of the same bundle, cf 
\cite{GelFuk}, \cite{Fuk}. 
Throughout this section we follow the notations of \cite{Fuk}. 

For any $k \geq 1$, let $M^m_k$ be the subset of $M^m$ of $m$-uples
$(x_1,\ldots,x_m)$ which do not have more then $k$ different 
entries $x_i$, that is, 
$$
	M^m_k = \{(x_1,\ldots,x_m)\in M^m\,\,|
	\,\,\forall(i_1,\ldots,i_{k+1})\subset(1,\ldots,m)\,\,
	\exists (i_r,i_s)\,\,:\,\,x_{i_r}=x_{i_s}\}. 
$$
Then, the submanifold $M^m_1=\{(x,\ldots,x)\in M^m \}=\triangle$ 
coincides with the diagonal of $M^m$, all the subsets $M^m_k$ for $k\geq m$ 
coincide with the whole $M^m$, and there is a sequence of inclusions 
$$
	0=M^m_0 \subset M^m_1 = \triangle \subset M^m_2 \subset \ldots 
	\subset M^m_{m-1} \subset M^m_m = M^m . 
$$ 

Denote by $CL^m_k(M)$ the subspace of $CL^m_{cont}(\Vect M)$ of 
generalized sections of the bundle $TM^{\otimes m}$ with support  
on the subset $M^m_k$. Then there is an induced sequence 
of inclusions of cochains modules  
$$
	\{0\}=CL^m_0(M) \subset CL^m_1(M) \subset\ldots\subset 
	CL^m_{m-1}(M) \subset CL^m_m(M)=CL^m_{cont}(\Vect M).  
$$ 
In other words, as a distribution, a Leibniz $m$-cochain is concentrated on 
some subspace $M^m_k$ (perhaps on the whole $M^m$).

The above inclusions of Leibniz cochain modules define an increasing 
multiplicative filtration of the complex $CL^*_{cont}(\Vect M)$ 
called {\em diagonal filtration\/}, that is,  
$$
	d(CL^m_k(M)) \subset CL^{m+1}_{k}(M), 
$$
for all $k$ and $m$, and  
$$
	CL^m_k(M) \cup CL^l_h(M) \subset CL^{m+l}_{k+h}(M), 
$$
for all $k,h$ and $m,l$.
To see this, it suffices to understand that the subspace $CL^m_k(M)$ 
consists exactely of Leibniz cochains which vanish on any family 
of $m$ vector fields having the property $(\triangle_k)$. 
(Recall that a family of $m$-vector fields $\Gamma \subset \Vect M$ 
has the property $(\triangle_k)$ if for any set of $k$ points $S \subset M$ 
at least one vector field of $\Gamma$ vanishes in a neighbourhood of $S$.)

Hence, there exists a spectral sequence abutting to the Leibniz cohomology 
of $\Vect M$. To avoid confusion, we denote this spectral sequence by 
$B_*(M)$. The $0$-th term is the bicomplex of the quotients 
$$
	B^{m,k}_0(M):= CL^{m}_{k}(M)/CL^{m}_{k-1}(M), 
$$
with differentials $d_{v}: B^{m,k}_0(M) \longrightarrow B^{m+1,k}_0(M)$ and 
$d_{t}: B^{m,k}_0(M) \longrightarrow B^{m+1,k-1}_0(M)$ (of bidegree $(1,0)$ 
and $(1,-1)$) induced by the Leibniz coboundary. 
For $k=1$, the differential $d_{v}$ coincides with the Leibniz
differential, hence the family $B^{m,1}_0(M)= CL^{m}_{\triangle}(M)$ 
is a subcomplex called {\em diagonal complex\/}. For $k=0$ we have 
$B^{m,0}_0(M)=0$, except for $m=0$ where $B^{0,0}_0(M)=\R$. 

Then the first term of the spectral sequence is the cohomology 
$$
\begin{array}{ll}
	B^{m,k}_1(M) = H^m(CL^{*}_k(M)/CL^{*}_{k-1}(M),d_v), 
	& \mbox{for $k>1$} \\ 
	B^{m,1}_1(M) = H^m(CL^{*}_{\triangle}(M),d_v), 
	& \mbox{for $k=1$} \\ 
	B^{m,0}_1(M) = \cases{\R, & $m=0$ \cr 0, & $m>0$ \cr}, 
	& \mbox{for $k=0$}. 
\end{array}
$$ 
The cohomology of the complex $(CL^{*}_{\triangle}(M),d_v)$ is also called 
the diagonal Leibniz cohomology of $\Vect M$, and denoted by 
$HL^*_{\triangle}(\Vect M)$. 
By abuse of notation, we call $k$-diagonal cochains the elements of 
$B^{*,k}_0(M)$, and $k$-diagonal cohomology the term $B^{*,k}_1(M)$, 
also denoted $HL^*_{(k)}(\Vect M)$. 
Of course, if $\tilde{d_{t}}$ is the differential induced by $d_{t}$ on 
the first term, the second term of the spectral sequence is the 
cohomology $B^{m,k}_2(M) = H^k(B^{m+*,*}_1(M),\tilde{d_t})$. 
The Leibniz cohomology of $\Vect M$ is then
$$
	HL^m(\Vect M) = \bigoplus_{k \geq m} B^{m,k}_{\infty}(M). 
$$

In this generality, as for the Gelfand-Fuks computations, the spectral
sequence is rather intractable. In order to compute its first term, we
introduce following Gelfand and Fuks, for each 
quotient complex $B^{*,k}_0(M)=CL^{*}_k(M)/C^{*}_{k-1}(M)$ 
a spectral sequence abutting to its cohomology 
$B^{m,k}_1(M)$. These spectral sequences are defined by the order 
filtration for Leibniz cochains, which is discussed in the next section. 


\section{The order filtration for diagonal cohomology}

The set $B^{m,k}_0(M)$ of $k$-diagonal Leibniz cochains contains the 
distributions concentrated on the subset $M^{m}_k$ modulo those 
concentrated on the subset $M^{m}_{k-1}$. 
The set $M^{m}_{(k)} = M^{m}_k \setminus M^{m}_{k-1}$ 
is now a submanifold without singularities, and a $k$-diagonal Leibniz 
cochain can then be described as a distributions defined on the jets of 
sections of the bundle $TM^{\otimes m}$, where the jet 
expansion is taken in the normal direction to the submanifold 
$M^{m}_{(k)}$ in $M^{m}$. 

Recall that a generalized section of a vector bundle $E$ over $M$, 
concentrated on a subset $S\subset M$, is of order $\leq l$ on $S$ 
if it vanishes on any section having trivial $l$-jet at every point of $S$. 

For any integer $p$, denote by $F^pB^{m,k}_0(M)$ the subspace of Leibniz 
$k$-diagonal $m$-cochains of $\Vect M$ which are of order $\leq m-p$. 
The family of such spaces defines a decreasing filtration of the complex 
$CL^{*}_k(M)/CL^{*}_{k-1}(M)$ of $k$-diagonal Leibniz cochains. 
In particular, for $k=1$, the family $F^pB^{m,1}_0(M)=F^pCL^m_{\triangle}(M)$ 
gives a filtration of the diagonal complex. 

\begin{theo}
\label{E=>B1}
For each $k>0$ there is a spectral sequence abutting to the $k$-diagonal 
Leibniz cohomology of $\Vect M$, with the following $E_2$-term:

\begin{displaymath}
	^{(k)}E_2^{p,q}=H_{-p}(M^k,M^k_{k-1})\otimes 
	\bigoplus_{\fila{q_1+\ldots+q_k=q}{q_1,\ldots,q_k>0}} 
	HL^{q_1}(W_n)\otimes\ldots\otimes HL^{q_k}(W_n), 
\end{displaymath}

\noindent 
where $H_{-p}(M^k,M^k_{k-1})$ denotes the relative cohomology of the 
manifold $M^k$ with respect to the subspace $M^k_{k-1}$. 

In particular, for $k=1$, there is a spectral sequence 
\begin{displaymath}
	E_2^{p,q} = H_{-p}(M)\otimes HL^q(W_n) \Rightarrow
	HL_{\triangle}^{p+q}(\Vect M) 
\end{displaymath}
\noindent 
abutting to the diagonal cohomology of $\Vect M$. 
\end{theo}

\noindent
{\bf Proof:}
Since the proof for Leibniz cohomology does not substantially differs from 
that for Lie cohomology, we only stress the main differences 
in the simple case of the diagonal complex. 

By definition, the $0$-th term of the spectral sequence defined by 
the order filtration on the diagonal complex $CL^*_{\triangle}(M)$ 
is the bicomplex 
$$
	E_0^{p,q}(M) = 
	F^pCL^{p+q}_{\triangle}(M)/F^{p+1}CL^{p+q}_{\triangle}(M)  
$$ 
that is, $E_0^{p,q}$ is the quotient space of cochains from
$CL^{p+q}_{\triangle}(\Vect M)$ which are of order $\leq q$ with
respect to the diagonal $\triangle$, modulo those of order $<q$. In other
words, an element of the space $E_0^{p,q}$ is a generalized section of
the bundle over $M$ 

\begin{displaymath}
{\cal E}^{p,q}_0={\rm Hom}(S^q{\rm norm}_{M^{p+q}}\triangle,
TM^{\otimes p+q}|_{\triangle}).
\end{displaymath}

\noindent
Here, ${\rm norm}_{M^{p+q}}\triangle$ is the normal bundle of the
diagonal (closed) submanifold $\triangle\subset M^{p+q}$, and is related 
to the fact that the jet expansion is in the normal direction with 
respect to the diagonal. The $q$th symmetric power $S^q$ arises from
the fact that the generalized sections are non-zero exactly on those
elements with non-trivial $q$-jet on the diagonal modulo those with
trivial $q+1$-jet. 
The fiber of the bundle ${\cal E}^{p,q}_0$ in a point $x \in M$ is

\begin{displaymath}
	{\rm Hom}(S^q(\bigoplus^{p+q}V/V_{\triangle}), V^{\otimes p+q}), 
\end{displaymath}

\noindent 
where $V=T_xM$ is the fiber of $TM$ and $V_{\triangle}$ denotes the image 
of the diagonal inclusion $V \hookrightarrow V\oplus\ldots\oplus V$.
The vector space $S^q(\bigoplus^{p+q}V/V_{\triangle})$ admits the 
standard Koszul resolution

\begin{displaymath}
\begin{array}{r} 
	0\gets S^q((\bigoplus^{p+q}V)/V_{\triangle})\stackrel{\ \ pr}{\gets}
	S^q(\bigoplus^{p+q}V)\gets S^{q-1}(\bigoplus^{p+q}V)\otimes
	V\gets\ldots \\ 
	\ldots\gets S^{q-i+1}(\bigoplus^{p+q}V)\otimes\Lambda^{i-1}V
	\gets\ldots 
\end{array}
\end{displaymath}

\noindent 
Using the isomorphisms 

\begin{displaymath}
\begin{array}{ll} {\displaystyle 
	{\rm Hom}(S^{q-i}(\bigoplus^{p+q}V)\otimes\Lambda^iV,V^{\otimes p+q})}  
	& {\displaystyle \cong \Lambda^iV'\otimes 
	{\rm Hom}(S^{q-i}(\bigoplus^{p+q}V),V^{\otimes p+q})} \\
	& {\displaystyle \cong \Lambda^iV'\otimes 
	{\rm Hom}(S^{q-i}(\bigoplus^{p+q}V) \otimes (V')^{\otimes p+q},\R), }
\end{array}
\end{displaymath}

\noindent 
and 
$$ 
	S^{*}(\bigoplus^{p+q}V) \otimes (V')^{\otimes p+q} 
	\cong (S^{*}(V) \otimes (V'))^{\otimes p+q} 
	\cong W_n^{\otimes p+q} ,  
$$
where $W_n$ is the Lie algebra of formal vector fields on $\R^n$, 
we finally have an isomorphism 
$$
	{\rm Hom}(S^{q-i}(\bigoplus^{p+q}V)\otimes\Lambda^iV,V^{\otimes p+q}) 
	\cong \Lambda^iV'\otimes CL^{p+q}_{(-p-i)}(W_n)
$$
where the index $l$ in $CL^{p+q}_{(l)}(W_n)$ means the weight with
respect to the adjoint action of the Euler field 
$e_{0}= \sum_{i=1}^n x_i \frac{d}{dx_{i}}$. 
Then the Koszul resolution gives rise to an exact sequence of fibres.
>From this, we get the corresponding exact sequence of bundles, and of sections,
because the global sections functor is exact on ${\cal C}^{\infty}$
fibre bundles. Finally, we pass to the exact sequence of generalized
sections by a lemma on duals of Fr\'echet nuclear spaces. 

In conclusion, we have an exact sequence of complexes (see \cite{Fuk}
p.146 for details)

\begin{displaymath}
0\gets E_0^{p,*}\gets Sec'\xi^{p,*}_0\gets Sec'\xi^{p,*}_1\gets\ldots
\end{displaymath}

\noindent
where $\xi^{p,q}_i$ the bundle associated with $TM$ with typical fiber
$\Lambda^i(TM)'\otimes\gamma_i^{p,q}$, and $\gamma_i^{p,q}$ 
is the bundle associated with $TM$ with typical fiber 
$CL_{(-p-i)}^{p+q}(W_n)'$. 
We use lemma (\ref{euler}) to get the
acyclicity of all complexes in the above sequence
except $E_0^{p,*}$ and $Sec'\xi^{p,*}_{-p}$.
The identification of the differential $d_1^{p,q}$ is clearly the same
as in \cite{GelFuk}, so we get the stated result. 
\hfill$\square$

This proof can easily be adapted to obtain the higher diagonal $E_2$-terms.


\section{Leibniz cohomology of the Lie algebra of vector fields on $S^1$}

Now we consider the circle $S^1$. We choose coordinates 
$z=\exp(i \vp)$, for $\vp \in [0,2\pi]$. 
As a topological Lie algebra, $\Vect S^1$ is spanned by the vector
fields $e_k :=i \exp(i k \vp)\frac{d}{d\vp}$, for all integers $k$, 
with Lie bracket $[e_k,e_l]=(k-l)e_{k+l}$. 

The Lie cohomology of $\Vect S^1$ was computed by Gelfand and Fuks. 
It is a free graded-commutative algebra with one even generator 
in degree $2$, represented by the Gelfand-Fuks cocycle $\omega$,  
and one odd generator in degree $3$, represented by the Godbillon-Vey 
cocycle $\theta$. Hence
$$
	H^{*}(\Vect S^1) = \bigoplus_{k \geq 0} \C[\omega^k] 
	\oplus \bigoplus_{k \geq 0} \C[\theta \cup \omega^k] .  
$$

For each $\vp_0 \in S^1$, the Taylor expansion of vector fields defines a map 
$$
\begin{array}{rcl}
	\pi_{\vp_0}: \Vect S^1 & \longrightarrow & W_1 \\ 
	f(\vp) \frac{d}{d\vp} & \longmapsto & 
	\sum_{n \geq 0} \frac{d^n f(\vp_0)}{d\vp^n}
	x^n \frac{d}{dx}, 
\end{array} 
$$
where $W_1 = \C[[x]] \otimes \{ \frac{d}{dx}\}$ is the complex Lie algebra 
of formal vector fields in the origin of $\R$ and $x=\vp - \vp_0$.  
The pull back $\pi^*$ gives a map from the Lie cocycles of $W_1$ to
the Lie cocycles of $\Vect S^1$, and similarly a map between the 
Leibniz cocycles. 
A cocycle $\gamma_{\vp_0} = \pi^*_{\vp_0} \gamma$ on $\Vect S^1$ which is 
pulled back from a cocycle $\gamma$ on $W_1$ is {\em local\/} on $S^1$, 
in the sense that as a distribution it has support on the single point 
$\vp_0$. 
Changing the point $\vp_0$ of evaluation produces cocycles which are 
cohomologous, hence a local cocycle admits an $S^1$-invariant form given 
by the integration 
$$
	\gamma := \int_0^{2\pi} \gamma_{\vp} d\vp, 
$$
where $S^1$ is meant to act on $\Vect S^1$ by rotations. 

The Godbillon-Vey cocycle $\theta$ on $\Vect S^1$ is precisely pulled
back from a cocycle on $W_1$, hence it is a local. 
Instead, the Gelfand-Fuks cocycle $\omega$ is a {\em diagonal\/} 
cocycle, that is, it is a distribution 
with support on the diagonal of $S^1 \times S^1$. 

The Leibniz cohomology of $W_1$, computed by Lodder, is the dual Leibniz 
algebra generated by the Godbillon-Vey cocycle $\theta$ and a new generator 
$\alpha$ in degree $4$. The cup products $\theta^2$ and $\alpha \cup \theta$ 
are zero, because of the definition of the cup product. 
The cup product $\alpha^2 \cup \alpha$ is equal to $2 \alpha \cup \alpha^2$, 
because of the relation of dual Leibniz algebra and the fact that 
$\alpha$ has even degree. 
Hence the only higher order Leibniz cocycles for $W_1$ are 
$\alpha^k := \alpha \cup \alpha^{k-1}$ and $\theta \cup \alpha^k$. 

The pull back on $\Vect S^1$ of these cocycles gives all the local
Leibniz cocycles which appear in $HL^{*}(\Vect S^1)$. We now determine
the diagonal ones.  

\subsection*{The diagonal Leibniz cohomology} 

The diagonal cohomology $HL^*_{\Delta}(\Vect S^1)$ is given by the cohomology 
classes which are represented by diagonal cocycles. 
In particular, it also contains local classes. 
The subset of local classes forms a dual Leibniz algebra with the cup product. 
In fact, the cup product of two local cocycles is still a local cocycle. 
Instead, the cup product of a local cocycle by a diagonal one, 
or the product of two diagonal cocycles, is not diagonal anymore, 
but $2$-diagonal. In general, the cup product of $k$ diagonal cocycle 
(with at least a non-local one) is $k$-diagonal. 
Therefore, the diagonal cohomology $HL^*_{\Delta}(\Vect S^1)$ is not a 
dual Leibniz algebra with the cup product, and it can only be described 
as a vector space.  

\begin{theo}
The diagonal cohomology of $\Vect S^1$ is the graded vector space 
spanned by the classes of the local cocycles 
$$
\begin{array}{cl} 
	\theta \cup \alpha^r & \mbox{in degree $3+4r$, for $r \geq 0$} \\ 
	\alpha^s & \mbox{in degree $4s$, for $s \geq 1$}
\end{array}
$$
and by the classes of the diagonal cocycles 
$$
\begin{array}{cl} 
	\omega_r & \mbox{in degree $2+4r$, for $r \geq 0$} \\ 
	\beta_s & \mbox{in degree $4s-1$, for $s \geq 1$}
\end{array}
$$
where $\omega_0 = \omega$ is the Gelfand-Fuks cocycle in degree $2$ 
and $\omega_{r}$, $\beta_{s}$ determine new invariants for $r,s \geq 1$. 
\end{theo}

\noindent
{\bf Proof:}
The spectral sequence of theorem (\ref{E=>B1}), which abuts to 
$HL^{*}_{\triangle}(\Vect S^1)$, degenerates at the second term 
for $M=S^1$, because 
$$
	E_2^{p,q}(S^1) = H_{-p}(S^1)\otimes HL^q(W_1)
$$ 
differs from $0$ only for $-p=0,1$, hence all the induced differentials 
in the higher terms are zero. 
So, we have 
$$
	HL^{m}_{\triangle}(\Vect S^1) = \bigoplus_{\fila{p+q=m}{q>0}} 
		H_{-p}(S^1)\otimes HL^q(W_1), 
$$ 
for all $m>0$, where the Leibniz cohomology of the formal vector
fields $W_1$ was computed by Lodder, in \cite{Lod}, as 
$$
	HL^q(W_1) = \bigoplus_{r+s=q}
	\Lambda^r[\theta]\otimes T^s[\alpha]. 
$$
Thus, we obtain 
$$
	HL^{m}_{\triangle}(\Vect S^1) = \bigoplus_{\fila{p+r+s=m}{r+s>0}} 
	\Lambda^p[\eta]\otimes \Lambda^r[\theta]\otimes T^s[\alpha], 
$$
where $\eta$ is a generator in degree $-1$, $\theta$ is the
Godbillon-Vey generator in degree $3$, and $\alpha$ is the Lodder
generator in degree $4$. 
Since $\Lambda^p[\eta]$ differs from $0$ only for $p=0,-1$, and 
$\Lambda^r[\theta]$ differs from $0$ only for $r=0,3$, we have 
\begin{displaymath}
\begin{array}{l} 
	HL^{m}_{\triangle}(\Vect S^1) = \displaystyle{ \bigoplus_{r+s=m+1}
	\R[\eta]\otimes \Lambda^r[\theta]\otimes T^s[\alpha] \oplus 
	\bigoplus_{r+s=m} \Lambda^r[\theta]\otimes T^s[\alpha]} \\ 
	\qquad = \displaystyle{\R[\eta \otimes \theta] \otimes
	T^{m-2}[\alpha] \oplus 
	\R[\theta] \otimes T^{m-3}[\alpha] \oplus 
	\R[\eta] \otimes T^{m+1}[\alpha] \oplus 
	T^{m}[\alpha]}.  
\end{array}
\end{displaymath}

\noindent 
If we call the new diagonal cocycles
$$
\begin{array}{ll}
	\omega_r:=\eta\otimes\theta\otimes\alpha^r, 
	&\mbox{for $r\geq 0$} \\ 
	\beta_s:=\eta\otimes\alpha^s, 
	&\mbox{for $s\geq 1$}, 
\end{array}
$$
and we remark that the higher-degree local cocycles are represented by 
the cup products, we get the final result 
$$
	HL^{*}_{\triangle}(\Vect S^1) = 
	\Lambda^*[\theta] \otimes T^*[\alpha] \oplus 
	\bigoplus_{r\geq 0} \C[\omega_{r}] \oplus 
	\bigoplus_{s\geq 1} \C[\beta_{s}].  
$$
\hfill$\square$

\begin{cor}
The diagonal Leibniz cohomology $HL^{*}_{\triangle}(\Vect S^1)$ is
periodic of period $4$. 
The dimensions of $HL^{n+4k}_{\triangle}(\Vect S^1)$ for $n=1,2,3,4$ are 
respectively $0,1,2,1$. 
\end{cor}

In particular, the second cohomology group $HL^{2}(\Vect S^1)$ has
dimension $1$, generated by the Gelfand-Fuks class. 
This result was obtained by Loday and Pirashvili in \cite{LP}.


\subsection*{Representative cocycles for the cohomology classes} 

Recall that the Godbillon-Vey class in $H^3(W_1)$ has a representative cocycle 
$$
	\theta(F,G,H) = \left| \begin{array}{lll} 
	f &g &h  \\ f' &g' &h' \\ f'' &g'' &h'' 
	\end{array} \right|_{x=0} ,    
$$
where $F=f(x)\frac{d}{dx}$, $G=g(x)\frac{d}{dx}$ and 
$H=h(x)\frac{d}{dx}$ denote some formal vector fields, and all 
the functions are evaluated at $x=0$. Now, there are two ways to
derive from $\theta$ the Gelfand-Fuks cocycle. 

Under the pull back $\pi^*_{\vp_0} : C^*(W_1) \longrightarrow C^*(\Vect S^1)$, 
$\theta$ determines a local cohomology class in $H^3(\Vect S^1)$. 
Its representative cocycle $\theta_{\vp_0}(F,G,H)$ can be expressed 
by the same local formula, where $F=f(\vp)\frac{d}{d\vp}$, $G$ and $H$ 
are now vector fields on $S^1$ and the functions are evaluated at $\vp = \vp_0$. 
In fact, the integral form of local cocycles has an advantage, since it 
may allow to determine automatically diagonal cocycles of one degree less, 
by contracting any of the variables with the Euler field $e_0 = i \frac{d}{d\vp}$. 
We explain this in detail. 
 
The only cochains which may contribute to the cohomology are those of zero weight, 
i.e. those $\gamma$ for which $ad_{e_0}(\gamma)=0$, because the subcomplex of 
non-zero weight cochains is contractible, as shown in appendix \ref{reduction}.
Among these, those of the form $\gamma= \gamma' \wedge \epsilon^0$, 
where $\epsilon^0 = e_0^*$ is the dual cochain of the Euler vector field and 
$\gamma'$ has zero weight, are such that the contraction gives 
$i_{e_0} (\gamma) = \gamma'$. 
Since the Lie coboundary operator satisfies Cartan's formula, we have 
$$
	d(\gamma) = d(\gamma') \wedge \epsilon^0 .  
$$ 
Hence, $\gamma$ is a cocycle if and only if $\gamma'$ is.  

However, if the local cocycle $\gamma$ is not in its integral form, the dependence 
$\gamma = \gamma_{\vp_0}$ from the evaluation point $\vp_0$ is reflected onto 
$\gamma' = \gamma'_{\vp_0}$, while the new cocycle $\gamma'$ is surely not local. 
Indeed, the term which leads to locality is precisely the $\epsilon^0$ that $\gamma'$, 
being antisymmetric, cannot contain. Since the contraction and the differential 
commute with the integration, we can avoid the dependence on the point of evaluation, 
by applying the procedure to the integral mean of $\gamma \wedge \epsilon^0$. 

The Godbillon-Vey cocycle is an example of a cocycle of the form
$\gamma \wedge \epsilon^0$, with $\gamma = -i/2 \ \epsilon^{-1} \wedge
\epsilon^1$ and $\epsilon^p = e_p^*$.  
For the Gelfand-Fuks class in $H^2(\Vect S^1)$, the well known representative cocycle 
$$
	\omega(F,G) = \int_{S^1} \left| \begin{array}{ll} 
	f'(\vp)&g'(\vp) \\ f''(\vp)&g''(\vp) \end{array} \right| d\vp    
$$
can be obtained from the Godbillon-Vey cocycle as follows, 
$$ 
	\omega = i_{e_0} \int_{S^1} \theta_{\vp} d\vp.  
$$ 
We apply the same method to get diagonal $k$-cocycles on $S^1$ from local 
Leibniz $k+1$-cocycles. 

The Godbillon-Vey class in $HL^3(W_1)$ has the same representative cocycle 
$\theta$, and by pull back we obtain the local class $\theta_{\vp_0}$ 
in $HL^3(\Vect S^1)$. 
By contracting its integral form we recover the Leibniz cocycle $\omega$ 
which represents the Gelfand-Fuks class.  

In fact, as a Leibniz cocycle, $\theta$ is also cohomologous to the 
cocycle 
$$
	\tilde{\theta}(F,G,H) = (f g h''' - f g''' h)|_{x=0}, 
$$
which, by integration and contraction, produces again the Gelfand-Fuks 
cocycle. However, we prefer to use the antisymmetric form of $\theta$ 
since it is advantageous in the computation of the cup products.

The Lodder class in $HL^4(W_1)$ is the image of the non trivial Lie 
cohomology class in $H^3(W_1;W'_1)$ with values in the adjoint representation, 
under the map $H^3(W_1;W'_1) \longrightarrow HL^3(W_1;W'_1) \cong H^4(W_1)$. 
It is represented by the cocycle 
$$
	\alpha(L,F,G,H) = l'(0) \left| \begin{array}{lll} 
	f &g &h \\ f' &g' &h' \\ f'' &g'' &h'' 
	\end{array} \right|_{x=0},  
$$
where $L = l(x)\frac{d}{dx}$ denote another formal vector field. 
Then, the first new class for the Leibniz cohomology of $\Vect S^1$ 
is represented by the cocycle 
$$
	\beta_1(L,F,G) = \int_{S^1} l'(\vp) \left| \begin{array}{ll} 
	f'(\vp) &g'(\vp) \\ f''(\vp) &g''(\vp) \end{array} \right| d\vp, 
$$
where, as before, $L, F, G$ now denote vector fields on $S^1$.

The second method to obtain the Gelfand-Fuks cocycle from the
Godbillon-Vey cocycle does not differ too much from the first one: one
can show that 

$$ \frac{d}{d\phi}\theta_{\phi} = d\left(\left|\begin{array}{cc} -' &
-' \\ -'' & -'' \end{array}\right|\right). $$ 

This permits to define $\omega$ as integral over $S^1$ of the
expression in paranthesis on the right hand side, because it gives
automatically a cocycle by the above equation. It is easily seen that
$\beta_1$ can also be obtained via this method, i.e. that in the
Leibniz setting 

$$ \frac{d}{d\phi}\alpha_{\phi} = d\left((-')\left|\begin{array}{cc} -' &
-' \\ -'' & -'' \end{array}\right|\right). $$  

The representative cocycles for the other new diagonal cohomology classes 
can be obtained in the same manner, starting from the cup products 
of the two local cocycles $\theta$ and $\alpha$. The formulas are though 
complicated by the presence of the shuffles, and we omit them. 

To end this section, let us just remark that the higher diagonal
Leibniz cohomology spectral sequences also collapse at the
second term. This is due to the fact that $H^*((S^1)^k,(S^1)^k_{k-1})$
is non-zero only in two dimensions, see \cite{GelFuk}.  


\section*{Aknowledgements}

The authors would like to thank Olivier Mathieu for many useful 
discussions and suggestions on the subject of the paper and Muriel 
Livernet for kindly remarking some oversights. 


\begin{appendix}
\section{Reduction to Euler-invariant cochains}
\label{reduction}

Suppose a Lie algebra $g$ possesses an element $e_0$, called the Euler field, 
and a basis of eigenvectors for the adjoint operator 
$ad_{e_0}: X \mapsto [X,e_0]$. 
For instance, let $\{ e_k, k \in \Z, k\neq 0 \}$ 
be the basis of elements such that $ad_{e_0}(e_k) = k e_k$. 

The adjoint operator $ad_{e_0}$ is a derivation of the Lie algebra $g$.  
It can be extended to the tensor powers of $g$ as a graded derivation 
(with respect to the tensor product), and consequently to the Leibniz 
cochains $CL^*(g)$. 

For any $\la \in \Z$, let $CL^*_{(\la)}(g)$ be the subset of Leibniz cochains 
$\gamma$ such that $ad_{e_0}(c) = \la c$.  
Formula $(iv)$ of proposition (3.1) in \cite{LP} means that the
adjoint operator $ad_{e_0}$ commutes with the Leibniz differential. 
Hence, $CL^*_{(\la)}(g)$ is a subcomplex for any $\la \in \Z$, and 
the Leibniz complex $CL^*(g) = \oplus_{\la} CL^*_{(\la)}(g)$ splits up 
into a direct sum of such subcomplexes. 
In particular, the elements $c$ of $CL^*_{(0)}(g)$ are called 
{\em Euler-invariant cochains\/}, because $ad_{e_0}(c) = 0$.

The same splitting, as a completed direct sum, occurs if the Lie algebra 
has a topology and the basis $\{ e_k, k \in \Z, k\neq 0 \}$ is a 
topological one. 

As for Lie cohomology, cf. theorem 1.5.2 \cite{Fuk}, we then have:  

\begin{lem}
\label{euler}
Suppose $g$ is a Lie algebra. Under the previous assumptions the 
Euler-invariant cohomology determines the Leibniz cohomology of $g$, 
that is, 

\begin{displaymath}
HL^*_{cont}(g) \cong H^*(CL^*_{(0)}(g)).
\end{displaymath}
\end{lem}

\noindent 
\pr We construct a contracting homotopy for the cochain complexes
$CL^*_{(\la)}(g)$  with $\la\not= 0$. 
For $\la\not= 0$, let $D^p_{(\la)}: CL^{p+1}_{(\la)}(g) \longrightarrow CL^p_{(\la)}(g)$ 
be defined by 

\begin{displaymath}
(D^p_{(\la)}c)(x_1,\ldots,x_{p-1}) = c(x_1,\ldots,x_{p-1},e_0)
\end{displaymath}

\noindent 
for any $p$-cochain $c$. Cartan's formula (formula $(i)$ in
proposition (3.1) \cite{LP}) shows that 

\begin{displaymath}
(d\tilde{D}^{p}_{(\la)}c)(x_1,\ldots,x_p) = (\la - \tilde{D}^{p+1}_{(\la)}d)(c)(x_1,\ldots,x_p).  
\end{displaymath}

\noindent where we have set $\lambda=\sum_{i=1}^p\lambda_i$ with
$ad_{e_0}(x_i)=:\lambda_ix_i$, and
$\tilde{D}^p_{(\lambda)}:=(-1)^pD^p_{(\lambda)}$. 
This means that $\tilde{D}^p_{(\la)}$ is a contracting homotopy for
$CL^*_{(\la)}(g)$ for $\la\not= 0$. \hfill$\square$
\bigskip
 
\end{appendix}


\end{document}